\newcommand{\bx}{{\bf x}}
\newcommand{\by}{{\bf y}}
\newcommand{\bz}{{\bf z}}

\newcommand{\ra}[1]{\renewcommand{\arraystretch}{#1}}

\documentclass{article}
\usepackage{newtxtext}
\usepackage[latin1]{inputenc}
\usepackage{amsmath}
\usepackage{amsfonts}
\usepackage{amssymb}
\usepackage{booktabs}
\usepackage{cancel}
\usepackage{cleveref}
\usepackage{colortbl}
\usepackage{csquotes}
\usepackage{float}
\usepackage{graphicx}
\usepackage{helvet}
\usepackage{latexsym}
\usepackage{listings}
\usepackage{mathtools}
\usepackage{mathpazo}
\usepackage{mathptmx}
\usepackage{mathrsfs}
\usepackage{multirow}
\usepackage{pgfplotstable}
\usepackage{siunitx}
\usepackage{tikz}
\usepackage{verbatim}
\usepackage{xcolor}
\usepackage{yfonts}
\begin{document}

\markboth{Naum Khutoryansky}
{Explicit Stencil Computation Schemes Generated by
	Poisson's Formula for the 2D Wave Equation.}

%
%

\title{EXPLICIT STENCIL COMPUTATION SCHEMES GENERATED BY
	POISSON'S FORMULA FOR THE 2D WAVE EQUATION}     

\author{Naum M. Khutoryansky \\ Dept. of Engineering Technology, Drexel University}
\date{}

\maketitle


\begin{abstract}

A new approach to building explicit time-marching stencil computation schemes for the transient
2D acoustic wave equation is implemented.
It is based on using Poisson's formula and its three time level modification combined with polynomial stencil interpolation of the solution at each time-step and exact integration.
The time-stepping algorithm consists of two explicit stencil computation procedures: a first time-step
procedure incorporating the initial conditions and a two-step scheme for the second and next time-steps.	
Three particular explicit stencil schemes (with five, nine and 13 space points) are constructed using this approach. Their stability regions are presented. Accuracy advantages of the new schemes in comparison with conventional finite-difference schemes are demonstrated by simulation using an exact benchmark solution.
\end{abstract}


	
\section{Introduction}

Stencil computations are widely implemented in many numerical algorithms that involve structured grids. In acoustic field simulation based on two- or three-dimensional transient wave equations, the finite difference time domain method is a standard approach leading to stencil operations \cite{Hamilton2015}. In an explicit time-stepping finite difference scheme, a solution value at each point in a time-space grid is calculated using a linear combination of values at its spatial neighbors from previous time steps. Among such schemes the main attention in the literature has been given to two-step schemes (which operate over three time levels $t_{k+1}=(k+1)\tau$, $t_{k}=k\tau$ and $t_{k-1}=(k-1)\tau$ where $\tau$ is a fixed time increment). They have been intensively studied and reviewed in many articles and books (see, e.g., $^{2-7}$). 

Explicit two-step numerical schemes for the scalar wave equation can also be devised based on spherical means representations as it is done in \cite{Alpert2000} where an integral evolution formula with three time levels is derived. For the 2D case, the evolution formula has been implemented in \cite{Li2004} and \cite{Hagstrom2015} using piecewise polynomial interpolation in 2D mesh cells and numerical integration. 

The two-step schemes pose some challenges when imposing the initial conditions. To calculate the value of a sought solution $u$ at the first time level $t_1$, one needs values of $u$ from time levels $t_0$ and $t_{-1}$. The initial condition for $u$ provides the required values at $t_0$. However, $u(\bx,t_{-1})$ must be inferred from the other initial condition. So, there is a need for transforming the two-step scheme at the first time step to a one-step form which has been discussed, e.g., in \cite{Strikwerda2004} and \cite{Langtangen2017}. A conventional approach for the 2D case (described, e.g., in \cite{Langtangen2017}) uses the central difference in time for approximating the initial condition for $v=\partial u / \partial t$.  For any point $\bx$ in $R^2$, the value of $u(\bx,t_{-1})$ is inferred as $u(\bx,t_1)-2 \tau v(\bx,0)$. Even though this approach is attractive due to its simplicity, it is worthwhile to consider the numerical schemes for both the first time step and the next steps derived using a unified approach. 

The present work focuses on building explicit time-stepping stencil computation schemes for the transient 2D acoustic (scalar) wave equation using spherical means formulas including Poisson's formula \cite{Evans1998} and a similar integral formula involving three time levels \cite{Alpert2000} that are described in section 2. A general form of the implemented time-stepping algorithm is presented in section 3 where two different integral expressions are given for the first time-step and for the next steps. Sections 4 and 5 consider some known results regarding polynomial interpolation on stencils and exact integration needed for deriving explicit stencil computation schemes from these integral expressions. Particular numerical schemes for five-, nine- and 13-point stencils are obtained in section 6. For each scheme, two separate expressions are derived: 1) a one-step expression for the first time-step through the stencil interpolated initial conditions; 2) a two-step expression for the second and next time-steps. All the obtained stencil expressions for the first time-step have not been previously presented in the literature. A derived two-step expression for the nine-point stencil is also new. In contrast, the obtained two-step expressions for five- and 13-point stencils coincide with the corresponding finite-difference stencil expressions. It is shown by simulation that the derived numerical schemes can significantly improve accuracy of stencil computations in comparison with conventional approaches.  

\section{Representation Formulas}
Consider a Cauchy problem for the transient 2D scalar wave equation
\begin{gather}  \label{we1}
\frac{\partial^2 u}{\partial t} -
c^2 \left( \frac{\partial^2 u}{\partial x_1^2}
+ \frac{\partial^2 u}{\partial x_2^2}
\right)=0,
\quad
u=u(\bx,t), \quad \bx=(x_1,x_2) \mbox{ in } \mathbb{R}^2,  
\\
u|_{t=0}=u_0(\bx), \quad \frac{\partial u}{\partial t}|_{t=0} =v_0(\bx).
\end{gather} 
Its solution is given by the representation formula which is often named Poisson's formula for the 2D wave  equation \cite{Evans1998}:
\begin{equation}   
u(\bx,t)=\frac{1}{2\pi c t^2} \int_{|\by-\bx|^2<t^2}
\frac{t u_0(\by)+t \ \nabla u_0(\by) \cdot (\by - \bx) +t^2 v_0(\by) }{\sqrt{c^2 t^2-|\by-\bx|^2}} d y_1 d y_2
\label{RF1}  
\end{equation}
where $\nabla$ denotes the gradient operator and $\by=(y_1,y_2) \in \mathbb{R}^2$ is a variable of integration.

Rewriting the previous integral on the unit disk, we get
\begin{equation}  \label{RF2} 
u(\bx,t)=\frac{1}{2\pi} \int_{|\bz|<1}
\frac{u_0(\bx+c t \bz)+c t \nabla u_0(\bx+c t \bz) \cdot \bz }{\sqrt{1-|\bz|^2}} d z_1 d z_2
+\frac{t}{2\pi} \int_{|\bz|<1}
\frac{v_0(\bx+c t \bz)}{\sqrt{1-|\bz|^2}} d z_1 d z_2
\end{equation}  
where $\bz=(z_1,z_2) \in \mathbb{R}^2$ is a new variable of integration.

Formulas \eqref{RF1} and \eqref{RF2} are also valid for negative $t$ which can be proved based on the time reversal property of the wave equation. 
By substituting $-t$ for $t$ in formula \eqref{RF2} and changing $\bz$ to $-\bz$ inside its integrals, we obtain the following analog of that formula for negative time:
\begin{equation}  \label{RF3}   
u(\bx,-t)=\frac{1}{2\pi} \int_{|\bz|<1}
\frac{u_0(\bx+c t \bz)+c t \nabla u_0(\bx+c t \bz) \cdot \bz }{\sqrt{1-|\bz|^2}} d z_1 d z_2  \\
-\frac{t}{2\pi} \int_{|\bz|<1}
\frac{v_0(\bx+c t \bz)}{\sqrt{1-|\bz|^2}} d z_1 d z_2.
\end{equation}
The only difference between the right-hand sides of \eqref{RF2} and \eqref{RF3} is the opposite signs of the second term. So, one can eliminate this term by summing \eqref{RF2} and \eqref{RF3}. Shifting the initial moment in the resulting formula from $t=0$ to $t=t_*$, we obtain the following expression involving three time points with a time increment $\tau$:
\begin{equation}
u(\bx,t_* + \tau)+u(\bx,t_* - \tau)=\frac{1}{\pi} \int_{|\bz|<1}
\frac{u(\bx+c t \bz,t_*)+c t \nabla u(\bx+c t \bz,t_*) \cdot \bz }{\sqrt{1-|\bz|^2}} d z_1 d z_2
\label{RF4}   
\end{equation}
where the right-hand side does not include the time derivative $v$. The same representation formula (in a different form) has been derived in \cite{Alpert2000} without using Poisson's formula.

Both formulas \eqref{RF2} and \eqref{RF4} will be used below to build a time-marching stencil computation algorithm for the wave equation.
\section{An Integral Time-Stepping Algorithm}

Consider a uniform time grid $\{t_0=0, t_1=\tau,\ldots, t_k=k\tau,\ldots\}$ where $\tau$ is a fixed time-step. Denote by $u_k(\bx)$ the restriction of $u(\bx,t)$ to a moment $t=t_k$. 
Next, denote by $A(\bx,\tau)$ and $B(\bx,\tau)$ the following integral operators acting on continuous functions defined in $\mathbb{R}^2$: 
\begin{equation}
A(\bx,\tau)f(\cdot)=\frac{1}{2 \pi} \int_{|\bz|<1}
\frac{f(\bx+c \tau \bz)+c \tau \nabla f(\bx+c \tau \bz) \cdot \bz }{\sqrt{1-|\bz|^2}} d z_1 d z_2,
\label{Op1}   
\end{equation}
\begin{equation}
B(\bx,\tau)f (\cdot)=\frac{\tau}{2 \pi} \int_{|\bz|<1}
\frac{f(\bx+c \tau \bz)}{\sqrt{1-|\bz|^2}} d z_1 d z_2
\label{Op2}   
\end{equation}

The time-stepping algorithm proposed here consists of two procedures based on the representation formulas \eqref{RF2} and \eqref{RF4}:

1) \textit{The procedure for the first time-step} which, according to \eqref{RF2}, calculates $u_1(\bx)$ as
\begin{equation}
u_1(\bx)=A(\bx,\tau)u_0(\cdot)+B(\bx,\tau)v_0(\cdot), \bx \in \mathbb{R}^2   
\label{A1}   
\end{equation}

2) \textit{The procedure for the second and next time-steps} which, using \eqref{RF4} for $t_*=t_k$, calculates $u_{k+1}(\bx)$  as
\begin{equation}
u_{k+1}(\bx)=2 A(\bx,\tau)u_k(\cdot)-u_{k-1}(\bx), \bx \in \mathbb{R}^2, k=1,2,... 
\label{A2}   
\end{equation}

While formula \eqref{A2} involves three time levels (two time-steps), formula \eqref{A1} for the first time-step incorporates only two time levels (one time-step) without using any finite difference approximation of the time derivative.

\section{Using Polynomial Interpolation on Stencils}

Consider a two-dimensional uniform Cartesian grid $\{(i_1h,i_2h)\}$ where $i_1$ and $i_2$ are integers, and $h$ is the grid spacing in both directions $x_1$ and $x_2$. Suppose that the evaluation point $\bx$ in formulas \eqref{A1} and \eqref{A2} is a grid point $\bx_{ij}=(ih,jh)$. Our intention is to choose a stencil in the Cartesian grid  and reduce the integral formulas \eqref{A1} and \eqref{A2} to linear combinations of the stencil node values of $u_k(\bx)\, (k=0,1,2,\ldots)$ and $v_0(\bx)$. Such a reduction will be done by using polynomial interpolation.

Assume that a particular stencil with $m$ nodes is chosen for polynomial interpolation. The corresponding index set $\{(q_1,q_2)\}$ is denoted by $Q_m$. The stencil index components $q_1$ and $q_2$ are numbered relative to the referencing point located at the evaluation point $\bx_{ij}$. So, polynomial interpolation in a neighborhood of the evaluation point $\bx_{ij}$ will be carried out using interpolation points
\begin{equation}
\bx_{i,j}+\bx_{q_1,q_2}=\bx_{i+q_1,j+q_2}, \quad (q_1,q_2) \in Q_m.
\end{equation}    

Following \cite{McKinney1972} we associate with $Q_m$ a set of $m$ distinct bivariate monomials
\begin{equation} \label{McKinney}  
\mathcal{M}_m= \{x_1^{\alpha(q_1)} x_2^{\alpha(q_2)}, (q_1,q_2) \in Q_m \}
\end{equation}
where 
\begin{equation}  \label{alpha} 
\alpha(q)= 
\begin{dcases}
2|q|-1  & \text{if } q < 0  \\
2q & \text{if } q\geq 0. 
\end{dcases}
\end{equation}
The above function has a unique inverse function
\begin{equation}  \label{invMcKinney} 
q(\alpha)=(-1)^\alpha \left[\frac{\alpha+1}{2}\right]
\end{equation}
where $[\cdot]$ is the whole part function.
So, for each index value $(q_1,q_2)) \in Q_m$ there exists a unique monomial from $\mathcal{M}_m$ and vice versa.

Consider a polynomial space $\mathcal{P}_m$ spanned by $\mathcal{M}_m$. We will use only those stencils for which the Lagrange interpolation problem is unisolvent in $\mathcal{P}_m$ (see \cite{Gasca2012}). In this case, there exists a Lagrange basis for $\mathcal{P}_m$ that can be built as described below. 

Suppose that there is an ordering imposed on the monomials in $\mathcal{M}_m$  
\begin{equation} \label{monomials}
\{\mu_1(\bx), \ldots, \mu_m(\bx)\}.
\end{equation}
The corresponding stencil nodes are numbered accordingly using \eqref{McKinney}-\eqref{invMcKinney}: 
\begin{equation}     \label{nodes}
\{\bx_{(1)}, ..., \bx_{(m)}\}.
\end{equation}
Thus, one can compute the following matrix  \cite{Gasca2012}:
\begin{equation} \label{matrix1}
D=\left[\mu_s(\bx_{(r)})\right]_{m \times m}
\end{equation}
If this matrix is non-singular, i.e., $\det(D) \ne 0$, which means that the Lagrange basis exists, then the inverse matrix 
\begin{equation} \label{matrix2}
C=\left[c_{sr}\right]_{m \times m}=D^{-1}
\end{equation}
can be calculated. Its components are instrumental in expressing the Lagrange basis functions through the chosen set of monomials:
\begin{equation} \label{Lagrange}
L_s(\bx)=\sum_{r=1}^m c_{sr} \mu_r (\bx), \quad s=1,...,m.
\end{equation}
Once the Lagrange basis is obtained, we again need two indexes to denote the Lagrange basis functions in accordance with the two index notation for grid points. The sequence \eqref{monomials} defines a relationship $s=g(\alpha_1,\alpha_2)$ between the ordinal number $s$ and monomial exponents $(\alpha_1,\alpha_2)$. Then, the resulting relationship between an index pair $(q_1,q_2)$ and the corresponding ordinal number $s$ is given by the following expression:
\begin{equation} \label{Lagrange3}
s=\gamma(q_1, q_2)=g(\alpha(q_1),\alpha(q_2)).
\end{equation}
where $\alpha(\cdot)$ is specified in \eqref{alpha}. 
Therefore, by introducing a new (two index) notation for the Lagrange basis functions
\begin{equation} \label{Lagrange4}
\phi_{q_1,q_2}(\bx)=L_{\gamma(q_1,q_2)}(\bx),
\end{equation}
we get the interpolation formula for a continuous function $f(\bx)$ in the form
\begin{equation} \label{Lagrange5}
\tilde{f}(\bx)=\sum_{(q_1, \, q_2) \in Q_m} f_{i+q_1,j+q_2} \phi_{q_1 q_2} (\bx) 
\end{equation}  

Even though different sets of monomials and their sequences can be employed for building Lagrange bases, we will use a particular method of monomial ordering that is described below. 

Denote by $\mathcal{M}$ the set of all monomials $x_1^{\alpha_1} x_2^{\alpha_2}$ where $\alpha_1$ and  $\alpha_2$  are natural numbers including $0$. For each monomial $x_1^{\alpha_1} x_2^{\alpha_2}$, the corresponding ordinal number $s$ will be assigned using the following function:
\begin{equation}  \label{g}
s=g(\alpha_1,\alpha_2 )=\frac{(\alpha_1+\alpha_2)(\alpha_1+\alpha_2+1)}{2}
+\begin{dcases}
\alpha_1 - \alpha_2 & \text{if } \alpha_2 < \alpha_1, \\
\alpha_2 - \alpha_1 + 1   & \text{if } \alpha_2 \geq \alpha_1. 
\end{dcases}
\end{equation}
It is easy to prove that the function \eqref{g} provides a one-to-one correspondence between $\mathcal{M}$ and the set of all positive natural numbers with the usual ordering. The order induced by \eqref{g} uses the total degree as the first sorting parameter (similarly to the more common graded lexicographic order) while the difference between individual degrees with the same total degree is used as the next sorting variable.

 The set $\mathcal{M}$ endowed with the order induced by \eqref{g} will be denoted by $\mathcal{M}^*$. The initial segment of $\mathcal{M}^*$ with $m$ members is denoted below by $\mathcal{M}_{\le m}^*$.  It will be shown in section 6 that such ordered monomial sets  play a useful role in building particular numerical schemes. 

\section{Calculating the Integrals}
Let $\bx=\bx_{ij}$ in the integral operators \eqref{Op1} and \eqref{Op2}. Without loss of generality, assume that the origin of the 2D coordinate system is located at $\bx_{ij}$ which can be achieved by a parallel translation of coordinates. Let $f(\by)$ appearing in \eqref{Op1} and \eqref{Op2}  be a monomial in local scaled variables $y_1/h$ and $y_2/h$:
\begin{equation} \label{monomial2}
\mu(\by)=\left(\frac{y_1}{h}\right)^{\alpha_1} \left(\frac{y_2}{h}\right)^{\alpha_2} , \quad \alpha_1,\alpha_2 = 0,1, \ldots.
\end{equation}
Then the integrals in \eqref{Op1} and \eqref{Op2}  can be exactly calculated and expressed through the Courant number
\begin{equation} \label{Courant}
\lambda=\frac{c\tau}{h}
\end{equation}
Indeed, using the table of integrals of Gradshteyn \& Ryzhik \cite{Gradshteyn1965}, expressions \eqref{Op1} and \eqref{Op2} are reduced to the following exact values:
\begin{equation} \label{exact1}
A(\bx_{ij},\tau)\mu(\cdot)=0, \, B(\bx_{ij},\tau) \mu(\cdot)=0 \, \mbox{  if $\alpha_1$ or $\alpha_2$ are non-negative odd integers},
\end{equation}
\begin{equation} \label{exact2}
\begin{aligned}
A(\bx_{ij},\tau) \mu(\cdot)=\frac{(\alpha_1-1)!! (\alpha_2-1)!!}{(\alpha_1+\alpha_2-1)!!} \lambda^{\alpha_1+\alpha_2}, \\ 
B(\bx_{ij},\tau) \mu(\cdot)=\frac{\tau}{\alpha_1+\alpha_2+1} A(\bx_{ij},\tau) \mu(\cdot) \\ 
\mbox{if $\alpha_1$ and $\alpha_2$ are both non-negative even integers},
\end{aligned}
\end{equation}
where $(\cdot)!!$ is the double factorial. It is assumed that $(-1)!!=1,0!!=1$.

The above formulas allow one to exactly calculate integrals \eqref{Op1} and \eqref{Op2} when $f(\cdot)$ is a polynomial from the Lagrange basis (see section 4).
\section{Particular Explicit Two-Step Schemes}
Now we can start building some numerical schemes by transforming the procedures \eqref{A1} and \eqref{A2} into algebraic expressions. All the functions $u_0(\bx), v_0(\bx)$ and $u_k(\bx), k=1,2,...$ included in these procedures will be interpolated in a stencil's center neighborhood using the same stencil nodes. The following standard notations for grid values of the solution and initial conditions will be used in the computation schemes:
\begin{equation} \label{grid0}
u_{ij}^0=u(\bx_{ij},0),v_{ij}^0 =v(\bx_{ij},0),
\end{equation}
\begin{equation} \label{grid1}
u_{ij}^k = u(\bx_{ij},k\tau), k=1,2,\ldots.
\end{equation}
\subsection{The five-point stencil}
Consider building the Lagrange basis for the space of complete second degree polynomials. In this case, the monomial basis sequence ordered according to \eqref{g} is as follows:
\begin{equation} \label{mono6}
\mathcal{M}_{\le 6}^*=\lbrace 1,x_1,x_2,x_1 x_2,x_1^2,x_2^2 \rbrace
\end{equation}
The corresponding interpolation stencil node sequence will be written according to \eqref{invMcKinney} as
\begin{equation*} \label{index6}
\lbrace(0,0),(-h,0),(0,-h),(-h,-h),(h,0),(0,h)\rbrace
\end{equation*}
The matrix $D$ for this stencil is non-singular with $\det(D)=4 h^8$. As a result, we get the Lagrange basis as
\begin{equation} \label{LB6}
\begin{split} 
\phi_{0,0}=L_1 &= 1+\frac{x_1 x_2}{h^2}-\frac{x_2^2}{h^2}-\frac{x_1^2}{h^2}, \quad
\phi_{-1,0}=L_2=-\frac{x_1}{2h}-\frac{x_1 x_2}{h^2}+\frac{x_1^2}{2h^2}, \\
\phi_{0,-1}&=L_3 =-\frac{x_2}{2h}-\frac{x_1 x_2}{h^2}+\frac{x_2^2}{2h^2},  \quad
\phi_{-1,-1}=L_4 = \frac{x_1 x_2}{h^2}, \\ 
&\phi_{1,0}=L_5 =\frac{x_1}{2 h}+\frac{x_1^2}{2 h^2}, \quad
\phi_{0,1}=L_6 =\frac{x_2}{2 h}+\frac{x_2^2}{2 h^2}.
\end{split}
\end{equation}
One can see that the Lagrange basis \eqref{LB6} includes monomials in the scaled variables $x_1/h$ and $x_2/h$ with coefficients independent of $h$. 

Using formulas from section 5 and notation \eqref{Lagrange4} we get
\begin{eqnarray} \label{A6}
A(\bx_{ij},\tau)\phi_{00}(\cdot)=1-2 \lambda ^2,  \quad A(\bx_{ij},\tau)\phi_{-1,-1}(\cdot)=0, \nonumber \\
A(\bx_{ij},\tau)\phi_{\pm 1,0}(\cdot)=A(\bx_{ij},\tau)\phi_{0, \pm 1}(\cdot)=\frac{1}{2} \lambda ^2,
\end{eqnarray}
\begin{eqnarray} \label{B6}
B(\bx_{ij},\tau)\phi_{0,0}(\cdot)=\tau (1-\frac{2}{3}\lambda^2),  \quad B(\bx_{ij},\tau)\phi_{-1,-1}(\cdot)=0, \nonumber \\
B(\bx_{ij},\tau)\phi_{\pm 1,0}(\cdot)=B(\bx_{ij},\tau)\phi_{0, \pm 1}(\cdot)=\frac{\tau}{6} \lambda ^2,
\end{eqnarray}
So, all coefficients for node $(-h,-h)$ disappear and the corresponding computational scheme contains only 5 spatial points which is shown in Figure 1.
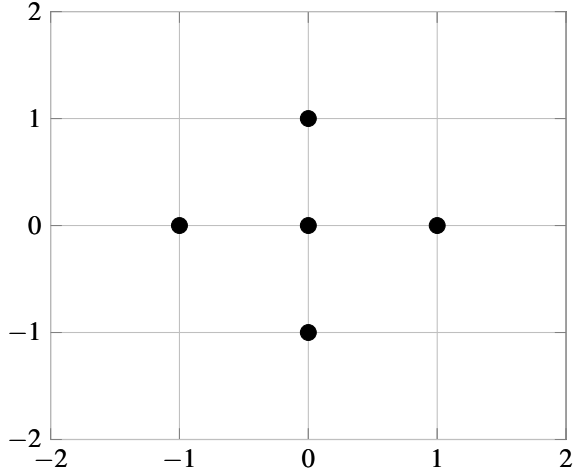
\begin{figure}[ht]
	\centering
	\begin{tikzpicture}
	\begin{axis}[xmin=-2,   xmax=2,
	ymin=-2,   ymax=2,
	extra x ticks={-2,-1,0,1,2},
	extra y ticks={-2,-1,0,1,2},
	extra tick style={grid=major}]
	\addplot [only marks,mark=*, mark size=3pt] 
	coordinates 
	{(0,0) (-1,0) (1,0) (0,-1) (0,1)};
	\end{axis}
	\end{tikzpicture}
	\caption{Index set for the five-point numerical scheme}
	\label{fig:Q5}
\end{figure}
\subsubsection{A new first time-step expression}
As a result, the proposed five-point numerical scheme is as follows:

1) \textit{for the first time-step}

\begin{multline} \label{comp5.1}
u_{ij}^1 = u_{ij}^0+\tau v_{ij}^0 
+\frac{\lambda^2}{2} \left( u_{i-1,j}^0+u_{i+1,j}^0+u_{i,j-1}^0+u_{i,j+1}^0-4 u_{ij}^0 \right) \\ 
+\frac{\tau\lambda^2}{6} \left(v_{i-1,j}^0+v_{i+1,j}^0+v_{i,j-1}^0
+v_{i,j+1}^0-4 v_{ij}^0 \right);
\end{multline} 
2) \textit{for the second and next time-steps}
\begin{equation} \label{comp5.2}
u_{ij}^{k+1} = 2 u_{ij}^k - u_{ij}^{k-1}
+\lambda^2 \left( u_{i-1,j}^k+u_{i+1,j}^k+u_{i,j-1}^k+u_{i,j+1}^k-4 u_{ij}^k \right) ,\quad k=1,2,\ldots.
\end{equation} 

A conventional approach for the five-point stencil (see, e.g., \cite{Langtangen2017}) uses the same procedure \eqref{comp5.2} for the second and next time-steps. However, for the first time-step,  this approach uses another formula (rather than \eqref{comp5.1}) based on the central difference in time for approximating the initial condition for $v$ combined with \eqref{comp5.2} for $k=0$:
\begin{equation} \label{comp5.1a}
u_{ij}^1 = u_{ij}^0+\tau v_{ij}^0
+\frac{\lambda^2}{2} \left( u_{i-1,j}^0+u_{i+1,j}^0+u_{i,j-1}^0+u_{i,j+1}^0-4 u_{ij}^0 \right).
\end{equation} 
Comparing \eqref{comp5.1} and \eqref{comp5.1a} one can see that the first three terms of the right-hand parts of these equations coincide, but the forth term present in the right-hand side of \eqref{comp5.1} is absent in \eqref{comp5.1a}. Thus, the difference between these two first time-step expressions depends on properties of $v_0(\bx)$. If, for example, $v_0(\bx)$ is a linear function of spatial coordinate, then there is no difference between  \eqref{comp5.1} and \eqref{comp5.1a}. However, for more general cases encountered in practice, the difference may exist and can be influential from the accuracy point of view which is shown in subsection 6.1.2.

The stability condition for numerical schemes utilizing \eqref{comp5.2} is well known from von Neumann stability analysis \cite{Strikwerda2004}:
\begin{equation} \label{stability5}
\lambda \le \lambda_{\max}=\frac{\sqrt 2}{2}.
\end{equation} 
\subsubsection{Numerical simulation}
To make a numerical comparison of both approaches, we use the initial and boundary conditions corresponding to a standing wave exact solution:
\begin{equation} \label{exact sol}
u_e(x_1,x_2,t)=\sin(2 \pi x_1) \sin(2 \pi x_2) \sin(2\sqrt2 \pi c t).
\end{equation} 
This solution creates the following pair of initial conditions for the numerical simulation:
\begin{equation} \label{Exact IC}
u_0(x_1,x_2)=0, \quad v_0(x_1,x_2)=2\sqrt2\pi c \ \sin(2 \pi x_1) \sin(2 \pi x_2).
\end{equation} 
We consider the unit square $\Omega=[0,1]^2$ as a space region for the numerical solution and apply boundary condition $u=0$ on $\partial\Omega$ generated by \eqref{exact sol}. In addition, we assume that $c=1$. Let $n_t$ be the number of time-steps and $n \geq 2$ be the space discretization number related to $h$ by equation $nh=1$. The above numerical schemes have been employed over a spatial index set $\{(i,j), 0<i,j<n\}$
using grid boundary conditions
\begin{equation*}
u_{0,j}^k=u_{n,j}^k=u_{i,0}^k=u_{i,n}^k=0, \quad 0 \leq i,j \leq n, 
\quad k \geq 0.
\end{equation*}
The accuracy of both numerical schemes have been estimated using the relative $L^2$ error defined as
\begin{equation} \label{error}
E(n,n_t)=\left(\frac{\sum_{k=1}^{n_t}\sum_{i=0}^{n}\sum_{j=0}^{n} \left(u_{ij}^k-u_e(ih,jh,k\tau)\right)^2}{\sum_{k=1}^{n_t}\sum_{i=0}^{n}\sum_{j=0}^{n} \left(u_e(ih,jh,k\tau)\right)^2}\right)^{1/2}
\end{equation}

The relative $L^2$ errors for the proposed new five-point scheme \eqref{comp5.1}$\&$\eqref{comp5.2} (based on Poisson's formula) and the conventional one \eqref{comp5.1a}$\&$\eqref{comp5.2} are denoted by $E_{P5}(n,n_t)$ and $E_{C5}(n,n_t)$, respectively. Calculated values of $E_{P5}(n,n_t)$, $E_{C5}(n,n_t)$ and their ratio for different combinations of $n$, $n_t$  and $\lambda=0.707$ are presented in table 1 below. For every value of $n$ used, three values of $n_t$ are considered: $n_t=1$ (to check errors after the first step), $n_t=n$ and $n_t=2n$ (to demonstrate the error accumulation process).
\begin{table}[h]
	\centering
	\ra{1.3}
	\caption{Wave \eqref{exact sol} simulations using the five-point stencil.}
	\begin{tabular}{@{}cccccc@{}}\toprule			
		$n \quad$ & $n_t \quad$ & $\lambda$ & $E_{P5} \quad$ & $E_{C5}$  \\   \midrule
		10 &  1 & $0.707$ & $9.0843\cdot10^{-4}$   & $6.8938\cdot10^{-2}$   \\
		10 & 10 & $0.707$ & $9.1540\cdot10^{-4}$   & $6.8945\cdot10^{-2}$   \\	
		10 & 20 & $0.707$ & $9.1604\cdot10^{-4}$   & $6.8945\cdot10^{-2}$   \\ \midrule	
		20 &  1 & $0.707$ & $5.4767\cdot10^{-5}$   & $1.6636\cdot10^{-2}$   \\
		20 & 20 & $0.707$ & $5.6800\cdot10^{-5}$   & $1.6638\cdot10^{-2}$   \\
		20 & 40 & $0.707$ & $5.7372\cdot10^{-5}$   & $1.6638\cdot10^{-2}$   \\ \midrule	
		40 &  1 & $0.707$ & $3.3924\cdot10^{-6}$   & $4.1230\cdot10^{-3}$   \\
		40 & 40 & $0.707$ & $4.0331\cdot10^{-6}$   & $4.1234\cdot10^{-3}$   \\
		40 & 80 & $0.707$ & $4.4928\cdot10^{-6}$   & $4.1234\cdot10^{-3}$   \\ \midrule
		80 &  1 & $0.707$ & $2.1158\cdot10^{-7}$   & $1.0285\cdot10^{-3}$   \\
		80 & 80 & $0.707$ & $4.3820\cdot10^{-7}$   & $1.0286\cdot10^{-3}$   \\
		80 &160 & $0.707$ & $6.5824\cdot10^{-7}$   & $1.0286\cdot10^{-3}$   \\ \bottomrule	
	\end{tabular}	
\end{table}
Table 1 demonstrated a much higher accuracy of the new scheme in comparison with the conventional one (with the error ratio $E_{C5}/E_{P5}$ exceeding $10^3$ for more dense grids).  Observed increases of the relative $L^2$ errors from $n_t=n$ to $n_t=2n$ are negligible (less than $10^{-6}$). 
\subsection{A nine-point square stencil} 
Consider the initial segment $\mathcal{M}_{\le 11}^*$ of $\mathcal{M}^*$ with the last member equal to $x_1^2 x_2^2$:
\begin{equation} \label{monomials11}
\{1,x_1,x_2,x_1 x_2,x_1^2,x_2^2,x_1^2 x_2,x_1 x_2^2,x_1^3, x_2^3,x_1^2 x_2^2\}
\end{equation}
This set includes all the monomials of total degree $\le 3$ and one monomial of total degree $=4$, and is the minimal initial segment of $\mathcal{M}^*$ that includes the bivariate tensor-product of the monomial bases for the second degree polynomials in each spatial coordinate. The corresponding index set $Q_{11}$ is shown in Figure 2 where both solid dots and empty circles denote interpolation points. 
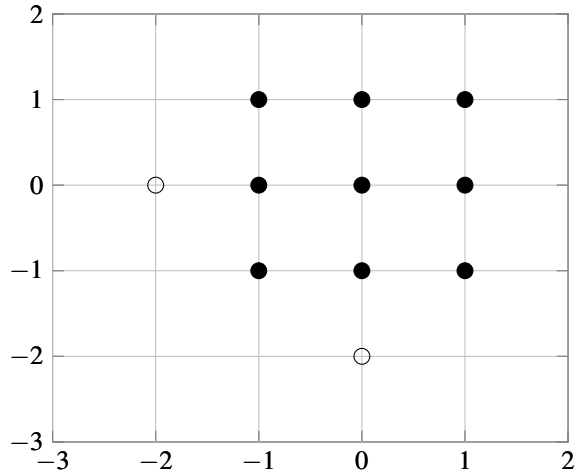
\begin{figure}[b!]
	\centering
	\begin{tikzpicture}
	\begin{axis}[xmin=-3,   xmax=2,
	ymin=-3,   ymax=2,
	extra x ticks={-3,-2,-1,0,1,2},
	extra y ticks={-3,-2,-1,0,1,2},
	extra tick style={grid=major}]
	\addplot [only marks,mark=*, mark size=3pt] 
	coordinates {(0,0) (-1,0) (0,-1) (-1,-1) (1,0) (0,1) (1,-1) (-1,1) (1,1)};
	\addplot [only marks,mark=o, mark size=3pt] 
	coordinates {(-2,0) (0,-2)};
	\end{axis}
	\end{tikzpicture}
	\caption{Index set for the 11-point interpolation stencil}
	\label{fig:Q9}
\end{figure}
The corresponding Lagrange basis can be easily calculated using matrix \eqref{matrix1} and presented similar to \eqref{LB6}. However, it is preferable to avoid presenting long expressions that include 11 monomials and display only terms that will be used later.  Denote by $\phi^e_{q_1,q_2}$ a part of $\phi_{q_1,q_2}$ that includes all the monomial terms with even exponents in both coordinates. Then we get
\begin{equation} \label{LB11}
\begin{split}
\phi^e_{0,0}=1-\frac{x_1^2}{h^2}-\frac{x_2^2}{h^2}+\frac{x_1^2 x_2^2}{h^4}, \quad &\phi^e_{-2,0}=\phi^e_{0,-2}=0, \\
\phi^e_{\pm 1,0}=\frac{1}{2}\left(\frac{x_1^2}{h^2}- \frac{x_1^2 x_2^2}{h^4} \right), \quad
\phi^e_{0,\pm 1}=\frac{1}{2}&\left(\frac{x_2^2}{h^2}- \frac{x_1^2 x_2^2}{h^4} \right), \quad \phi^e_{\pm 1,\pm 1}=\frac{1}{4}\frac{x_1^2 x_2^2}{h^4}.   
\end{split}
\end{equation}
\subsubsection{A new nine-point time-stepping scheme}
It follows from the results \eqref{exact1} and \eqref{exact2} of section 5 that functions $\phi^e_{q_1,q_2}$ rather than the complete Lagrange basis will be used in building numerical schemes. Therefore, two points (-2,0) and (0.-2) (presented by empty circles in Figure 2) will disappear in the corresponding numerical scheme. The remaining nodes (solid circles in Figure 2) create the nine-point square-shaped computational stencil. Using \eqref{exact1}-\eqref{exact2} we obtain:
\begin{equation} \label{A11}
\begin{split}
A(\bx_{ij},\tau)\phi_{00}(\cdot)&=1-2 \lambda ^2+\frac{\lambda ^4}{3},  \quad A(\bx_{ij},\tau)\phi_{\pm 1,\pm 1}(\cdot)= \frac{\lambda ^4}{12}, \\
A(\bx_{ij},\tau)&\phi_{\pm 1,0}(\cdot)= A(\bx_{ij},\tau)\phi_{0,\pm 1}(\cdot)=\frac{\lambda ^2}{2}-\frac{\lambda ^4}{6},
\end{split}
\end{equation}
\begin{equation} \label{B11}
\begin{split}
B(\bx_{ij},\tau)\phi_{00}(\cdot)&=\tau\left(1-\frac{2\lambda ^2}{3}+\frac{\lambda ^4}{15}\right),  \quad B(\bx_{ij},\tau)\phi_{\pm 1,\pm 1}(\cdot)= \frac{\tau\lambda ^4}{60}, \\
B(\bx_{ij},\tau)&\phi_{\pm 1,0}(\cdot)= B(\bx_{ij},\tau)\phi_{0,\pm 1}(\cdot)=\frac{\tau\lambda ^2}{6} \left(1-\frac{\lambda ^2}{5}\right).
\end{split}
\end{equation}
Now we can use \eqref{A1} and \eqref{A2} to build a new nine-point numerical scheme similar to  \eqref{comp5.1}-\eqref{comp5.2}. However, to avoid writing long expressions, some additional notations will be needed:
\begin{equation} \label{notations2}
\begin{split}
\delta^k_{ij}(q_1,q_2)=u^k_{i+q_1,j+q_2}+u^k_{i-q_2,j+q_1}
+u^k_{i1-q_1,j-q_2}
& +u^k_{i+q_2,j-q_1}-4 u^k_{i,j},  \\
\epsilon^0_{ij}(q_1,q_2)=v^0_{i+q_1,j+q_2}+v^0_{i-q_2,j+q_1}
+v^0_{i-q_1,j-q_2}
& +v^0_{i+q_2,j-q_1}-4 v^0_{i,j}, \\ 
&k=0,1,2,\ldots.
\end{split}
\end{equation}
Thus, the following time-stepping numerical scheme is derived:

1) \textit{for the first time-step}
\begin{equation} \label{comp9.1}
\begin{split}
u_{ij}^1 = u_{ij}^0+\tau v_{ij}^0+\frac{\lambda^2}{2} \Bigg[ \left(1-\frac{\lambda^2}{3}\right) \delta^0_{ij}(1,0)+\frac{\lambda^2}{6}\delta^0_{ij}(1,1) \Bigg] \\ 
+\frac{\tau\lambda^2}{6}\Bigg[\left(1-\frac{\lambda^2}{5}\right) \epsilon^0_{ij}(1,0)+\frac{\lambda^2}{10}\epsilon^0_{ij}(1,1)
\Bigg];
\end{split}
\end{equation} 

2) \textit{for the second and next time-steps}
\begin{multline} \label{comp9.2}
u_{ij}^{k+1} = 2 u_{ij}^k - u_{ij}^{k-1}+\lambda^2 \Bigg[ \left(1-\frac{\lambda^2}{3}\right) \delta^k_{ij}(1,0)+\frac{\lambda^2}{6}\delta^k_{ij}(1,1)\Bigg],\quad k=1,2, \ldots.
\end{multline}
Using the von Neumann stability analysis method (see, e.g., \cite{Strikwerda2004}), we obtain the stability condition for the above numerical scheme as 
\begin{equation} \label{stability9}
\lambda \le \lambda_{\max}=\frac{\sqrt {3-\sqrt 3}}{2} \approx 0.796.
\end{equation} 

To our knowledge this numerical scheme has not been presented previously in the literature. A conventional explicit nine-point square-shaped scheme (dubbed the isotropic scheme (see, e.g., \cite{Trefethen1982}) has a different form obtained using a nine-point finite difference approximation of the Laplace operator in the two-dimensional space \cite{Kantorovich1958}:
\begin{equation} \label{comp9.2b}
u_{ij}^{k+1} = 2 u_{ij}^k - u_{ij}^{k-1}+\lambda^2\Bigg[\frac{2}{3} \delta^k_{ij}(1,0) +\frac{1}{6}\delta^k_{ij}(1,1)\Bigg],\quad k=1,2, \ldots.
\end{equation}
	
Comparing the coefficients in \eqref{comp9.2} and \eqref{comp9.2b} it is easy to see that the term $\delta^k_{ij}(1,1)$ is less influential in the new nine-point scheme than in the conventional nine-point scheme (with a coefficient ratio equal to $\lambda^2$). The stability condition for the conventional scheme is less restrictive than for the new one: $\lambda \le \lambda_{\max}=\sqrt {3}/{2} \approx 0.866$. On the other hand, the new scheme has some accuracy advantages in comparison with the conventional nine-point approach for $\lambda \le 0.796$ which is shown in the next subsection.
\subsubsection{Numerical simulation}
We consider simulation results for the proposed nine-point scheme \eqref{comp9.1}-\eqref{comp9.2} using initial and boundary conditions generated by the exact solutions of the previous subsection. Let $c=1$. We assume that the unit square $\Omega=[0,1]^2$ is used as a space grid region, the time grid interval is $[0,\lambda]$ with $n_t=n$. A comparison is made with simulated results on the same grids for the conventional scheme \eqref{comp9.2b}. 
Even though no expression for the first time-step corresponding to the conventional nine-point scheme is presented in \cite{Trefethen1982},  the usual approach based on the central difference for approximating the initial condition for $v$ combined with \eqref{comp9.2b} provides the corresponding expression
\begin{equation} \label{comp9.1b}
u_{ij}^1 = u_{ij}^0+\tau v_{ij}^0 +\frac{\lambda^2}{2}\Bigg[ \frac{2}{3}\delta^0_{ij}(1,0) +\frac{1}{6} \delta^0_{ij}(1,1) \Bigg] 
\end{equation} 
which is used for numerical simulation. 

The simulation results are presented in the next table for standing wave simulations. We have used two values of $\lambda$ for the simulation: $\lambda=0.707$, as in subsection 6.1, and $\lambda=0.796$, according to the stability condition \eqref{stability9}. Relative $L^2$ errors for the schemes \eqref{comp9.1}$\&$\eqref{comp9.2} and  \eqref{comp9.1b}$\&$\eqref{comp9.2b} are denoted by $E_{P9}$ and $E_{C9}$, respectively. According to the simulated data, the new scheme \eqref{comp9.1}$\&$\eqref{comp9.2} appears to be more accurate than the other scheme for both $\lambda=0.707$ and $\lambda=0.796$. 
\begin{table}[h]
	\centering
	\ra{1.3}
	\caption{Wave \eqref{exact sol} simulations using nine-point stencils}
	\begin{tabular}{@{}ccccc@{}}\toprule			
		$ n $ & $ n_t $ & $\lambda$ & $E_{P9}$  & $E_{C9}$ \\   \midrule          		
		10 & 10 & $0.707$    & $3.7058\cdot10^{-2}$ 
		 &$1.1741\cdot10^{-1}$ \\	
		10 & 10 & $0.796$    & $2.9587\cdot10^{-2}$ &  $1.1241\cdot10^{-1}$ \\	 \midrule  	
		20 & 20 & $0.707$    & $8.9333\cdot10^{-3}$ &  $2.8002\cdot10^{-2}$ \\
		20 & 20 & $0.796$    & $8.0697\cdot10^{-3}$ &  $2.7523\cdot10^{-2}$ \\  \midrule  	
		40 & 40 & $0.707$    & $2.3723\cdot10^{-3}$ &  $6.8821\cdot10^{-3}$ \\
		40 & 40 & $0.796$    & $2.5737\cdot10^{-3}$ &  $6.8668\cdot10^{-3}$ \\   \midrule  
		80 & 80 & $0.707$    & $7.5573\cdot10^{-4}$ &  $1.7084\cdot10^{-3}$ \\	
		80 & 80 & $0.796$    & $1.0274\cdot10^{-3}$ &  $1.7187\cdot10^{-3}$ \\ \bottomrule	
	\end{tabular}	
\end{table}

 The obtained results show a loss in accuracy of the both nine-point schemes for the simulated problem in comparison with the new five-point scheme results presented in Table 1.
\subsection{A 13-point stencil}
The next numerical stencil is based on complete bivariate interpolation polynomials of the fourth degree. Consider the initial segment $\mathcal{M}_{\le 15}^*$ of $\mathcal{M}^*$:
\begin{equation} \label{monomials15}
\{1,x_1,x_2,x_1 x_2,x_1^2,x_2^2,x_1^2 x_2,x_1 x_2^2,x_1^3, x_2^3,x_1^2 x_2^2, x_1^3 x_2, x_1 x_2^3, x_1^4,x_2^4 \}
\end{equation}
This set includes all the monomials of total degree $\le 4$.
The corresponding index set $Q_{15}$ is presented in Figure 3 where both thirteen solid dots and two empty circles denote interpolation points. 
\begin{figure}[ht]
	\centering
	\begin{tikzpicture}
	\begin{axis}[xmin=-3,   xmax=3,
	ymin=-3,   ymax=3  ,
	extra x ticks={-3,-2,-1,0,1,2,3},
	extra y ticks={-3,-2,-1,0,1,2,3},
	extra tick style={grid=major}]
	\addplot [only marks,mark=*, mark size=3pt] 
	coordinates {(0,0) (-1,0) (0,-1) (-1,-1) (1,0) (0,1) (1,-1) (-1,1) (1,1) (-2,0) (2,0) (0,-2) (0,2)};
	\addplot [only marks,mark=o, mark size=3pt] 
	coordinates {(-2,-1) (-1,-2)};
	\end{axis}
	\end{tikzpicture}
	\caption{Index set for the 15-point interpolation stencil}
	\label{fig:Q13}
\end{figure}
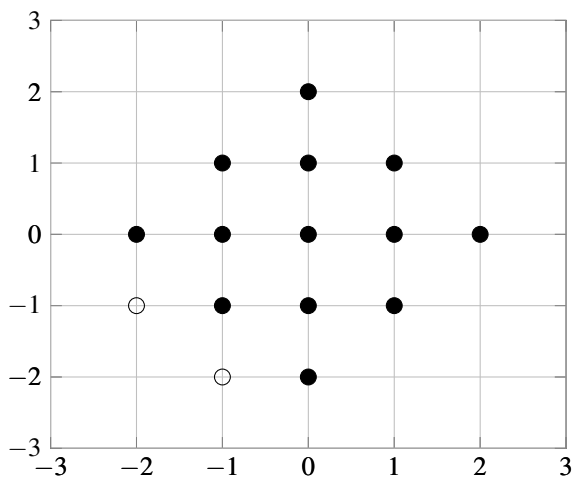

After calculating the corresponding Lagrange basis (see section 4) and using \eqref{exact1}-\eqref{exact2}, one can determine that all the coefficients in expressions \eqref{A1}-\eqref{A2} related to the two empty circles will disappear. As a result, we get a numerical scheme that involves the 13-point stencil (the solid dots in Figure 3). The notations \eqref{notations2} will be used in presenting the scheme to make the expressions more compact:

1) \textit{for the first time-step}
\begin{equation} \label{comp13.1}
\begin{split}
u_{ij}^1 = u_{ij}^0+\tau v_{ij}^0+\frac{\lambda^2}{2} \Bigg[ \frac{4-2\lambda^2}{3} \delta^0_{ij}(1,0)
+\frac{\lambda^2}{6}\delta^0_{ij}(1,1) 
+\frac{\lambda^2-1}{12}\delta^0_{ij}(2,0) \Bigg] \\ 
+\frac{\tau\lambda^2}{6}\Bigg[ \left(\frac{4}{3}-\frac{2\lambda^2}{5}\right) \epsilon^0_{ij}(1,0)
+\frac{\lambda^2}{10}\epsilon^0_{ij}(1,1)
+\left(\frac{\lambda^2}{20}-\frac{1}{12}\right)\epsilon^0_{ij}(2,0)
\Bigg];
\end{split}
\end{equation} 

2) \textit{for the second and next time-steps}
\begin{equation} \label{comp13.2}
\begin{split}
u_{ij}^{k+1} = 2 u_{ij}^k - u_{ij}^{k-1}
+\lambda^2 \Bigg[\frac{4-2\lambda^2}{3} \delta^0_{ij}(1,0)
+\frac{\lambda^2}{6}\delta^0_{ij}(1,1) 
+\frac{\lambda^2-1}{12}\delta^0_{ij}(2,0) \Bigg],  \\
k=1,2, \ldots.
\end{split}
\end{equation} 
Formula \eqref{comp13.2} for the second and next time-steps completely coincides with that obtained previously \cite{Cohen1987,Cohen1996} using the finite-difference method. 

However, expression \eqref{comp13.1} for the first time-step  has not been presented in the literature before. An advantage of using this expression rather than the conventional one (based on the central difference for approximating the initial condition for $v$) is shown in Table 3 where simulation results for the standing wave \eqref{exact sol} solution are presented. The relative $L^2$ errors for the conventional and new (Poisson's) approaches are denoted by $E_{C13}$ and $E_{P13}$, respectively.  For simplicity sake, the periodic boundary conditions have been incorporated in the simulation taking into account that the solution \eqref{exact sol} is periodic in both spatial directions.
Since the maximal Courant number needed for stability of this scheme is $1/\sqrt{2}$, a value of $\lambda=0.707$ was used. 
\begin{table}[h!]
	\centering
	\ra{1.3}
	\caption{Wave \eqref{exact sol} simulations using the 13-point stencil.}
	\begin{tabular}{@{}cccccc@{}}\toprule			
		$n \quad$ & $n_t \quad$ & $\lambda$ & $E_{P13} \quad$ & $E_{C13}$   \\   \midrule
		10 & 10 & $0.707$ & $4.2146\cdot10^{-5}$   & $6.8938\cdot10^{-2}$   \\	\midrule	
		20 & 20 & $0.707$ & $6.6004\cdot10^{-7}$   & $1.6636\cdot10^{-2}$   \\  \midrule	
		40 & 40 & $0.707$ & $1.1471\cdot10^{-8}$   & $4.1230\cdot10^{-3}$   \\ \midrule
		80 & 80 & $0.707$ & $2.8884\cdot10^{-10}$   & $1.0285\cdot10^{-3}$   \\ \bottomrule	
	\end{tabular}	
\end{table}\textit{}

A large advantage in accuracy for the new scheme demonstrated in Table 3 can be attributed to a higher accuracy of the new first time-step expression \eqref{comp13.1} in comparison with the conventional one. On the other hand,  it is worth to notice that the relative $L^2$ error $E_{C13}$ (corresponding to the conventional first time-step approach) in Table 3 has almost the same values as the error $E_{C5}$ in Table 1 despite using a higher degree interpolation stencil in Table 3. That is, an error introduced at the first time-step probably suppresses advantages of using a higher degree interpolation at later time-steps.

\section{Summary}
A new method is implemented to build explicit time-stepping stencil computation schemes for the transient 2D acoustic wave equation. It is based on using Poisson's formula and a similar three time level expression combined with polynomial stencil interpolation of the solution at each time-step and exact integration. As a result, for each chosen 2D stencil and a set of monomials, a unified time-marching scheme is created that includes two explicit computation procedures: for the first time-step and for the next steps.     

Particular explicit stencil computation schemes (with five, nine and 13 space points) are derived. All of the obtained first time-step computation expressions are different from those used in conventional finite-difference methods. The obtained two-step stencil expressions for the five- and 13-point stencils (where the complete interpolation polynomials have been used)  coincide with the corresponding finite difference schemes. The obtained two-step expression for the nine-point stencil is new. Its stability region is determined by the von Neumann analysis.

Simulation comparison results are presented for a benchmark problem with an exact solution. It is demonstrated by simulation that the proposed stencil computation approach maintains an accuracy advantage in comparison with conventional finite difference schemes which is mostly attributed to the new first time-step computation expressions.

\end{document}